\documentclass[12pt]{article}
\usepackage[cp1251]{inputenc}
\usepackage[russian]{babel}
\usepackage{amsfonts}
\usepackage{amssymb,amsmath,amsthm,eufrak,euscript}
\usepackage{graphicx,fancybox}

\newtheorem{theorem}{Теорема}

\begin{document}

\begin{center}
{\Large On the integrable magnetic geodesic flow on a 2-torus}
\end{center}

\medskip

\begin{center}
{\large S.V. Agapov}
\footnote{The work was supported by RSF (grant 14-11-00441).}
\end{center}

\begin{quote}
\noindent{\sc Abstract. } In this paper the magnetic geodesic flow on a 2-torus is considered. We study a semi-hamiltonian quasi-linear PDEs which is equivalent to the existence of polynomial in momenta first integral of magnetic geodesic flow on fixed energy level. It is known that diagonal metric associated with this system is Egorov one if degree of the first integral is equal to 2 or 3. In this paper we prove this fact in the case of existence of the first integral of any degree. \medskip

 \end{quote}

\section{\large Introduction}

We will consider a magnetic geodesic flow on a 2-torus.
We fix an energy level and assume that there is an additional first integral which is polynomial in momenta. It is known that the corresponding quasi-linear PDEs is semi-hamiltonian (see [1]), i.e. in the hyperbolic region it possesses Riemann's invariants and can be presented in the form of conservation laws. For any semi-hamiltonian system there is a diagonal metric naturally associated with this system. In [1] it is proved that if there is an integral of the second or of the third degree, then this metric is the one of Egorov type. In this work we generalize this result on the case of an arbitrary degree.

Firstly let me remind some known results about geodesic flow in absence of magnetic field. There exist two types of metrics on a 2-torus with an integrable geodesic flow. If metric is of the kind
$ds^2 = \Lambda (\alpha x + \beta y) (dx^2 + dy^2)$ or
$ds^2 = (\Lambda_1 (\alpha_1 x + \beta_1 y) + \Lambda_2 (\alpha_2 x + \beta_2 y)) (dx^2 + dy^2),$
then there exists an additional first integral which is polynomial in momenta of the first or of the second degree. It is not known if there exist metrics with irreducible polynomial first integrals of higher degree. This question was studied in [2] -- [5].
If the geodesic flow is integrable, then one can introduce global semi-geodesic coordinates $(t,x)$ (see [6]) on a 2-torus such that
$$
ds^2 = g^2(t,x)dt^2+dx^2, \qquad H = \frac{1}{2} \left( \frac{p_1^2}{g^2}+p_2^2 \right).
$$
The first integral has the form
$$
F = \frac{a_0}{g^n} p_1^n + \frac{a_1}{g^{n-1}} p_1^{n-1}p_2 + \ldots + \frac{a_{n-2}}{g^2} p_1^2 p_2^{n-2} + \frac{a_{n-1}}{g}p_1p_2^{n-1} + a_np_2^n, \ a_k = a_k(t,x).
$$

The condition $\dot{F} = \{F,H\} = 0$ is equivalent to the quasi-linear PDEs of the following form
$$
U_t+A(U)U_x=0 \eqno(1)
$$
on the coefficients of $F$. Here $U = (a_0, \ldots, a_{n-2}, a_{n-1})^T, \ a_{n-1}=g, \ a_n=1,$
\[ A = \left( \begin{array}{ccccccccc}
0 & 0 & \ldots & 0 & 0 & a_1\\
a_{n-1} & 0 & \ldots & 0 & 0 & 2a_2-na_0\\
0 & a_{n-1} & \ldots & 0 & 0 & 3a_3-(n-1)a_1\\
\ldots & \ldots & \ldots & \ldots & \ldots & \ldots\\
0 & 0 & \ldots & a_{n-1} & 0 & (n-1)a_{n-1}-3a_{n-3}\\
0 & 0 & \ldots & 0 & a_{n-1} & na_n-2a_{n-2}
\end{array} \right). \]
The system (1) is semi-hamiltonian (see [6]). It means that it can be written in the form of conservation laws, i.e. there exists a change of variables
$
U^T \rightarrow (G_1(U), \ldots, G_n(U))
$
such that for some $F_1(U), \ldots, F_n(U)$ the following relations hold:
$$
(G_j(U))_t + (F_j(U))_x = 0, \qquad j =1, \ldots, n.
$$
Moreover, in the hyperbolic region where all the eigenvalues $\lambda_1, \ldots, \lambda_n$ of matrix $A$ are real and pairwise distinct the system (1) possesses Riemann's invariants, i.e. there exists a change of variables
$$
U^T \rightarrow (r_1(U), \ldots, r_n(U))
$$
such that (1) can be written in the following form
$$
(r_j)_t + \lambda_j(r)(r_j)_x = 0, \qquad j =1, \ldots, n.
$$
Semi-hamiltonian systems were introduced and studied by S.P. Tsarev in [7], [8] (see also [9]).

The question of existence of an additional polynomial in momenta first integral of an arbitrary degree of geodesic flow in conformal coordinates $ds^2 = \Lambda (x,y) (dx^2+dy^2)$ was studied in [10]. The existence of the first integral of the form
$$
F = a_0 p_1^n + a_1 p_1^{n-1} p_2 + \ldots +a_n p_2^n, \qquad a_k = a_k(x,y)
$$
(herewith due to Kolokoltsov's theorem (see [5]) the following relations hold true
$$
a_n = c_1+a_{n-2}-a_{n-4}+ \ldots, \qquad a_{n-1} = c_2+a_{n-3}-a_{n-5}+ \ldots,
$$
where $c_1, c_2$ are some constants) leads to a quasi-linear PDEs of the kind
$$
A(U)U_x+B(U)U_y=0, \eqno(2)
$$
where $U = (a_0, \ldots, a_{n-2}, \Lambda)^T$. Systems of such a kind were studied, for example, in [11]. The system (2) is also semi-hamiltonian (in the regions where at least one of the matrixes $A$ and $B$ is nondegenerate).

Let me remind that for any semi-hamiltonian system the following relations on the eigenvalues hold true:
$$
\partial_{r_j} \frac{\partial_{r_i} \lambda_k}{\lambda_i-\lambda_k} = \partial_{r_i} \frac{\partial_{r_j} \lambda_k}{\lambda_j-\lambda_k}, \qquad i \neq j \neq k \neq i.
$$
It means that there exists a diagonal metric
$$
ds^2 = H_1^2 (r) dr_1^2 + \ldots + H_n^2 (r) dr_n^2 \eqno(3)
$$
such that its Christoffel symbols satisfy the following relations:
$$
\Gamma_{ki}^k = \frac{\partial_{r_i} \lambda_k}{\lambda_i-\lambda_k}, \qquad i \neq k.
$$

It is proved in [10] that the metric (3) associated with the system (2) is the one of Egorov type, i.e. its rotation coefficients $\beta_{kl}$ are symmetric:
$$
\beta_{kl} = \beta_{lk}, \qquad \beta_{kl} = \frac{\partial_{r_k} H_l}{H_k}, \qquad k \neq l,
$$
or, equivalently, there exists a function $a(r)$ such that
$\partial_{r_k} a(r) = H_k^2(r).$
Here $H_i$ are Lame coefficients of the metric (3), $H_i^2 = g_{ii}.$
Following [12], we shall call the corresponding semi-hamiltonian systems the Egorov ones.

According to the Pavlov--Tsarev theorem (see [12]), if the system is not breaking up $(\partial_{r_i} \lambda_k \neq 0, \ i \neq k)$, then it is Egorov one iff it possesses two conservation laws of the special form:
$$
F_x+G_y=0, \qquad F_y+H_x=0.
$$
In [10] these conservation laws are found explicitly for the system (2).

In this work we obtain the analogous results for the magnetic geodesic flow.

\section{\large The main theorem}

Consider the Hamiltonian system
$$
\dot{x}^j = \{x^j,H\}_{mg}, \qquad \dot{p}_j = \{p_j,H\}_{mg}, \qquad j=1,2 \eqno(4)
$$
on a 2-torus in magnetic field with Hamiltonian
$H = \frac{1}{2} g^{ij} p_ip_j$
and the Poisson bracket of the following kind:
$$
\{F,H\}_{mg} = \sum_{i=1}^2 \left ( \frac{\partial F}{\partial x^i} \frac{\partial H}{\partial p_i} - \frac{\partial F}{\partial p_i} \frac{\partial H}{\partial x^i} \right ) + \Omega (x^1,x^2) \left ( \frac{\partial F}{\partial p_1} \frac{\partial H}{\partial p_2} - \frac{\partial F}{\partial p_2} \frac{\partial H}{\partial p_1} \right ).
$$
If $\{F,H\}_{mg} =0$, then the function $F$ is the first integral of the geodesic flow (4). Magnetic geodesic flows (or, equivalently, systems with gyroscopic forces) were studied, for example, in [13] -- [16].

Choose the conformal coordinates $(x,y)$ in which
$ds^2 = \Lambda (x,y) (dx^2+dy^2)$, $H = \frac{p_1^2+p_2^2}{2 \Lambda}.$ Fix the energy level $H = \frac{1}{2}.$ Then one can parameterize the momenta by the following way:
$$
p_1 = \sqrt{\Lambda} \cos \varphi, \qquad p_2 = \sqrt{\Lambda} \sin \varphi.
$$
The equations (4) take the form
$$
\dot{x} = \frac{\cos \varphi}{\sqrt{\Lambda}}, \qquad \dot{y} = \frac{\sin \varphi}{\sqrt{\Lambda}}, \qquad \dot{\varphi} = \frac{\Lambda_y}{2 \Lambda \sqrt{\Lambda}} \cos \varphi - \frac{\Lambda_x}{2 \Lambda \sqrt{\Lambda}} \sin \varphi - \frac{\Omega}{\Lambda}.
$$

Following [1], we shall search for the first integral $F$ of the kind
$$
F(x,y,\varphi) = \sum_{k=-N}^{k=N} a_k(x,y) e^{i k \varphi}. \eqno(5)
$$
Here $a_k= u_k+ i v_k, a_{-k} = \bar{a}_k.$
The condition $\dot{F} = 0$ is equivalent to the following equation
$$
F_x \cos \varphi + F_y \sin \varphi + F_{\varphi} \left( \frac{\Lambda_y}{2 \Lambda} \cos \varphi - \frac{\Lambda_x}{2 \Lambda} \sin \varphi - \frac{\Omega}{\sqrt{\Lambda}} \right) = 0. \eqno(6)
$$

Let's substitute (5) into (6) and equate the coefficients at $e^{i k \varphi}$ to zero. We obtain
$$
\frac{\Lambda_y}{2 \Lambda} \frac{i (k-1) a_{k-1} + i (k+1) a_{k+1}}{2} - \frac{\Lambda_x}{2 \Lambda} \frac{i (k-1) a_{k-1} - i (k+1) a_{k+1}}{2 i} +
$$
$$
+\frac{(a_{k-1})_x+(a_{k+1})_x}{2} + \frac{(a_{k-1})_y-(a_{k+1})_y}{2 i} - \frac{i k \Omega a_k}{\sqrt{\Lambda}} = 0, \eqno(7)
$$
where $k = 0, \ldots, N+1, \ a_k = 0$ while $k > N.$

After eliminating the magnetic field $\Omega$ (see below) we obtain a quasi-linear PDEs on $a_j$ and $\Lambda$ of the kind
$$
A(U)U_x+B(U)U_y=0, \eqno(8)
$$
where $U = (\Lambda, u_0, \ldots, u_{n-1}, v_1, \ldots, v_{n-1})^T.$ We shall not write it down here explicitly in view of its bulkiness. It's proved in [1] that (8) is the semi-hamiltonian system for any $N.$  It's also proved in [1] that in the case of $N=2,3$ the system (8) is the Egorov one. In this work we generalize this result on the case of an arbitrary $N.$

\begin{theorem}
The system (8) is the one of Egorov type for any $N.$
\end{theorem}


\section{The proof of the theorem 1}

To prove the theorem 1 we will need only some of the equations of the system (7). In the case of $k = N+1$ we obtain the relation
$$
(a_N \Lambda^{-\frac{N}{2}})_x - i (a_N \Lambda^{-\frac{N}{2}})_y = 0,
$$
in what follows that one can put $a_N = \Lambda^{\frac{N}{2}}$ (see [1]).

Put $k=N$ in (7) and consider the real and imaginary parts of the obtained equation. We obtain the following expression for the magnetic field:
$$
\Omega = \frac{(N-1) (\Lambda_y u_{N-1}-\Lambda_x v_{N-1}) + 2\Lambda ((v_{N-1})_x-(u_{N-1})_y)}{4N\Lambda^{\frac{N+1}{2}}}, \eqno(9)
$$
as well as the following relation:
$$
2 \Lambda ((u_{N-1})_x + (v_{N-1})_y) = (N-1) (v_{N-1} \Lambda_y + u_{N-1} \Lambda_x). \eqno(10)
$$

Introduce the new variables
$$
f_k = u_k\Lambda^{-\frac{k}{2}}, \qquad g_k = v_k\Lambda^{-\frac{k}{2}}, \qquad k = 0, \ldots, N-1.
$$
Then the relations on $a_k$ become simpler. It follows from (9), (10) that
$$
\Omega = \frac{(g_{N-1})_x-(f_{N-1})_y}{2 N}, \eqno(11)
$$
$$
(f_{N-1})_x+(g_{N-1})_y = 0. \eqno(12)
$$

Put $k = N-1$ in (7), one obtains the equations:
$$
(N-1) f_{N-1} ((g_{N-1})_x-(f_{N-1})_y) + N ((f_{N-2})_y-(g_{N-2})_x-N \Lambda_y) = 0,
$$
$$
(N-1) g_{N-1} ((g_{N-1})_x-(f_{N-1})_y) + N ((f_{N-2})_x+(g_{N-2})_y+N \Lambda_x) = 0,
$$
which can be written, due to (12), in the following form:
$$
R_x + \left( \frac{N-1}{2} (g_{N-1}^2-f_{N-1}^2) - N^2 \Lambda + N f_{N-2} \right)_y = 0, \eqno(13)
$$
$$
R_y + \left( \frac{N-1}{2} (f_{N-1}^2-g_{N-1}^2) - N^2 \Lambda - N f_{N-2} \right)_x = 0, \eqno(14)
$$
where $$R = (N-1) f_{N-1} g_{N-1} - N g_{N-2}.$$ In what follows that (8) is the Egorov system.


\medskip

\begin{center}
{\bf \large The bibliography}
\end{center}

[1] M. Bialy, A.E. Mironov,
\emph {New semi-hamiltonian hierarchy related to integrable magnetic flows on surfaces}, Central European Journal of Mathema\-tics, 10:5 (2012), 1596 -- 1604.

[2] M.V. Pavlov, S.P. Tsarev,
\emph {On Local Description of Two-Dimensional Geodesic Flows with a Polynomial First Integral}, arXiv: 1509.03084v1.

[3] V.V. Kozlov et al.,
\emph {Polynomial integrals of geodesic flows on a two-dimensional torus},Math. USSR Sb., 83:2 (1995), 469 -- 481.

[4] V.V. Kozlov, D.V. Treshchev,
\emph {On the integrability of Hamiltonian systems with toral position space}, Math. USSR Sb., 63:1 (1989), 121 -- 139.

[5] V.N. Kolokoltsov,
\emph {Geodesic flows on two-dimensional manifolds with an additional first integral that is polynomial in the velocities}, Math. USSR-Izvestiya., 21:2 (1983), 291 -- 306.

[6] M. Bialy, A. E. Mironov,
\emph {Rich quasi-linear system for integrable geodesic flow on 2-torus}, Discrete and Continuous Dynamical Systems - Series A, 29:1 (2011), 81 -- 90.

[7] S.P. Tsarev,
\emph {On Poisson brackets and one-dimensional hamiltonian systems of hydrodynamic type}, Doklady Mathematics, 31 (1985), 488 -- 491.

[8] S.P. Tsarev,
\emph {The geometry of Hamiltonian systems of hydrodynamic type. The generalized hodograph method.}, Math. USSR-Izvestiya, 37:2 (1991), 397 -- 419.

[9] B. Sevennec,
\emph {Geometrie des systemes de lois de conservation}, V. 56, Memories, Soc.Math.de France, Marseille, 1994.

[10] M. Bialy, A.E. Mironov,
\emph {Integrable geodesic flows on 2-torus: formal solutions and variational principle}, Journal of Geometry and Physics, 87:1 (2015), 39 -- 47.

[11] M. Bialy,
\emph {Richness or semi-hamiltonicity of quasi-linear systems that are not in evolution form}, Quarterly of Applied Math. 2013. V. 71, 787 -- 796.

[12] M.V. Pavlov, S.P. Tsarev,
\emph {Tri-Hamiltonian Structures of Egorov Systems of Hydrodynamic Type}, Func. Anal. and Its Appl., 37:1 (2003), 32 -- 45.

[13] V.V. Kozlov,
\emph {Symmetries, topology, and resonances in Hamiltonian mechanics}, Springer, Verlag, Berlin. 1996.

[14] V.V. Ten,
\emph {Polynomial first integrals for systems with gyroscopic forces}, Math. Notes, 68:1 (2000), 135 -- 138.

[15] S.V. Bolotin,
\emph {First integrals of systems with gyroscopic forces}, Vestnik Moskov. Univ. Ser. 1 Mat. Mekh., 6 (1984), 75 -- 82.

[16] I.A. Taimanov,
\emph {On an integrable magnetic geodesic flow on the two-torus}, arXiv: 1508.03745v1.

\medskip
\medskip
\medskip
\medskip
\medskip

S.V. Agapov

Sobolev Institute of Mathematics, Novosibirsk Russia

Novosibirsk State University, Novosibirsk Russia

{\it E-mail address:} \ {\bf agapov.sergey.v@gmail.com}

\end{document}